\documentclass[12pt]{amsart}
\usepackage[utf8]{inputenc}
\usepackage{tikz}
\usetikzlibrary{arrows.meta, calc}
\usepackage{graphicx}
\usepackage{float}
\usepackage[T1]{fontenc}
\usepackage{lmodern}
\usepackage{amsmath,amssymb,amsthm,mathtools}
\usepackage[a4paper,margin=27mm]{geometry}
\usepackage{microtype}
\usepackage{enumitem}
\usepackage[colorlinks=true,linkcolor=blue,citecolor=blue,urlcolor=blue]{hyperref}
\usepackage{tikz-cd}
\hypersetup{pdftitle={Holomorphic sectional curvature of flag manifolds via Kähler reduction},pdfauthor={Ching-Hsien Lee}}
\numberwithin{equation}{section}
\newtheorem{theorem}{Theorem}[section]
\newtheorem{proposition}[theorem]{Proposition}
\newtheorem{lemma}[theorem]{Lemma}

\theoremstyle{definition}

\allowdisplaybreaks[1]
\newcommand{\C}{\mathbb C}

\newcommand{\ii}{\mathrm{i}}
\newcommand{\g}{\mathfrak g}
\newcommand{\h}{\mathfrak h}
\newcommand{\m}{\mathfrak m}
\newcommand{\II}{\mathrm{II}}
\newcommand{\HS}{\mathrm{HS}}

\DeclareMathOperator{\Ad}{Ad}
\DeclareMathOperator{\ad}{ad}
\DeclareMathOperator{\Tr}{Tr}
\DeclareMathOperator{\Hom}{Hom}

\DeclareMathOperator{\diag}{diag}
\title[Holomorphic sectional curvature of flag manifolds]{Holomorphic sectional curvature of flag manifolds via K\"ahler reduction}
\author{Ching-Hsien Lee}
\address{Department of Mathematics, Rutgers University, Piscataway, NJ 08854, U.S.A.}
\email{cl1539@math.rutgers.edu}
\date{\today}
\keywords{Generalized flag manifold, holomorphic sectional curvature, K\"ahler reduction, totally geodesic embedding}
\begin{document}
\begin{abstract}
We give an alternative proof of Li and Petersen's theorem that every invariant K\"ahler metric on a compact generalized flag manifold has positive holomorphic sectional curvature. For a unitary partial flag manifold, we realize the given metric by K\"ahler reduction of a Hermitian vector space. The defining equations place the level set on a sphere of radius $\varrho$, where the radial component of its second fundamental form gives the explicit lower bound $4/\varrho^2$ through the quotient curvature formula. Finally, we use a holomorphic isometric totally geodesic embedding to transfer the bound to generalized flag manifolds of every Lie type.
\end{abstract}
\maketitle

\section{Introduction}
Let $M = G/H$ be a compact generalized flag manifold, where $G$ is a compact connected semisimple Lie group and $H$ is the centralizer of a torus in $G$. Fix a $G$-invariant complex structure $J$ on $M$. These are the homogeneous complex manifolds known as K\"ahler $C$-spaces~\cite[Theorems~1 and~2]{Borel}. Ni and Zheng~\cite[Conjecture~1.8]{NZ} recorded the folklore conjecture that the canonical K\"ahler--Einstein metric on every K\"ahler $C$-space has positive holomorphic sectional curvature. Itoh~\cite[Theorem~3.1]{Itoh} proved positivity when the second Betti number is one and the holomorphic isotropy representation has at most two irreducible summands. Lohove~\cite[Theorem~22]{Lohove} proved that every invariant K\"ahler metric has positive holomorphic sectional curvature when the isotropy representation has at most three irreducible summands. He also proved positivity for the canonical K\"ahler--Einstein metric on every flag manifold of type $G_2$~\cite[Theorem~47]{Lohove}. Li and Petersen~\cite[Theorem~1.1, Proposition~3.2 and Corollary~3.3]{LP} settled the conjecture by proving positivity for every invariant K\"ahler metric. Their proof expresses the numerator of the holomorphic sectional curvature as an integral of squared norms and gives an explicit positive lower bound.

Projective embeddings give another approach to curvature on flag manifolds. Takeuchi~\cite{Takeuchi} described homogeneous K\"ahler submanifolds of complex projective space as orbits of highest-weight lines. Nakagawa and Takagi~\cite[Section~7]{NT} studied projective K\"ahler submanifolds with parallel second fundamental form. Console and Fino~\cite[Lemma~3.1]{CF} expressed the second fundamental form of a homogeneous projective orbit in terms of the unitary representation. For the incidence flag of a line in a hyperplane in $\C^{n+1}$ with $n \geq 2$, Yang and Zheng~\cite[Proposition~4.1]{YZ} obtained the sharp curvature range $[1/2,2]$ for the metric induced from two Fubini--Study factors of holomorphic sectional curvature $2$.

Our proof is motivated by Lohove's use of totally geodesic embeddings to transfer positivity between classical K\"ahler--Einstein flags in the same Lie family~\cite[Theorems~52 and~53]{Lohove}. Loi, Mossa and Zuddas~\cite[Theorem~1.1]{LMZ} constructed holomorphic isometric embeddings of flag manifolds into classical flag manifolds. For the given invariant metric on $M$, we use a holomorphic isometric embedding into a unitary partial flag manifold and prove that its image is totally geodesic. The Gauss equation then identifies the holomorphic sectional curvature of $M$ with that of the unitary flag in the corresponding tangent directions.

To estimate the curvature of the unitary flag, we realize its metric by K\"ahler reduction of a Hermitian vector space $V$. Let $Z$ be the level set of the moment map and let $K$ be the group acting on $V$. Kobayashi's formula~\cite[Theorem~(1.3)]{Kobayashi} expresses the holomorphic sectional curvature of $Z/K$ in terms of the second fundamental form of $Z$ in $V$. Manikandan~\cite{Manikandan} computed the holomorphic sectional curvature of smooth K\"ahler quotients and obtained nonnegativity for moduli spaces of stable quiver representations. For unitary linear actions, Zhang~\cite[Theorem~3.1 and Lemma~3.3]{Zhang} expressed the quotient curvature in terms of the moment map and proved strict positivity when the level set is compact. In our construction, taking traces in the defining equations of $Z$ shows that every point of $Z$ has the same norm $\varrho$, which we compute from the eigenvalues and their multiplicities in the chosen representation. The radial component of the second fundamental form then gives the explicit lower bound $4/\varrho^2$.

Let $g$ be a $G$-invariant K\"ahler metric on $(M,J)$ and let $\omega(X,Y) = g(JX,Y)$ be its K\"ahler form. For the Levi--Civita connection $\nabla$ of $g$ and real vector fields $X_1,X_2,X_3$, our curvature convention is
\begin{equation}\label{eq:Riemannian-curvature-tensor}
    R(X_1,X_2)X_3 = \nabla_{X_1}\nabla_{X_2}X_3 - \nabla_{X_2}\nabla_{X_1}X_3 - \nabla_{[X_1,X_2]}X_3.
\end{equation}
For $p \in M$ and $0 \ne X \in T_pM$, the holomorphic sectional curvature is
\begin{equation}\label{eq:conventions}
    H_g(X) = \frac{g(R(X,JX)JX,X)}{g(X,X)^2}.
\end{equation}
Fix a faithful finite-dimensional unitary representation of $G$ on a Hermitian vector space $W$. Identify $G$ with its image in $U(W)$ and its Lie algebra $\g$ with the corresponding subalgebra of $\mathfrak u(W)$. On $\mathfrak u(W)$ use the inner product $Q(X,Y) = -\Tr(XY)$. Let $p_G\colon G\to G/H$ be the quotient map and let $e$ be the identity of $G$. For this choice of $Q$, the K\"ahler form $\omega$ determines a unique element $\xi_\omega \in \g$ such that
\begin{equation}\label{eq:representative}
    (p_G^*\omega)_e(X,Y) = -Q(\xi_\omega,[X,Y])
\end{equation}
for all $X,Y \in \g$, where $[X,Y]$ denotes the Lie bracket. The centralizer of $\xi_\omega$ in $U(W)$ is
\[
    C_{U(W)}(\xi_\omega) = \{u \in U(W) : u\xi_\omega u^{-1} = \xi_\omega\},
\]
and its centralizer in $G$ is $C_G(\xi_\omega) = G \cap C_{U(W)}(\xi_\omega)$. Set $M_W = U(W)/C_{U(W)}(\xi_\omega)$. Let $p_W\colon U(W)\to M_W$ be the quotient map and let $I_W$ denote the identity on $W$. Under the conjugation action of $U(W)$ on $\mathfrak u(W)$, the map $uC_{U(W)}(\xi_\omega) \mapsto u\xi_\omega u^{-1}$ identifies $M_W$ with the adjoint orbit of $\xi_\omega$. The invariant inner product $Q$ identifies this orbit with the coadjoint orbit through $Q(\xi_\omega,\cdot)$. Let $\omega_W$ be the pullback to $M_W$ of the Kirillov--Kostant--Souriau form on that coadjoint orbit. With our sign convention,
\begin{equation}\label{eq:ambient-representative}
    (p_W^*\omega_W)_{I_W}(U,V) = -Q(\xi_\omega,[U,V])
\end{equation}
for all $U,V \in \mathfrak u(W)$. By Mason and Ziegler~\cite[Section~3, Theorem~1]{MZ}, $\omega_W$ admits a compatible $U(W)$-invariant complex structure $J_W$. Denote the associated K\"ahler metric by $g_W = \omega_W(\cdot,J_W\cdot)$.

\begin{theorem}\label{thm:main}
The inclusion $G \subset U(W)$ induces a holomorphic, isometric, totally geodesic embedding $\Phi\colon (M,g) \to (M_W,g_W)$ given by
\[
    \Phi(aH) = aC_{U(W)}(\xi_\omega).
\]
The stabilizer of $p_W(I_W)$ in $G$ is $C_G(\xi_\omega) = H$. Suppose that $\dim_{\C}M > 0$. Let $a_1 > \cdots > a_{r+1}$ be the distinct eigenvalues of the Hermitian endomorphism $-\ii\xi_\omega$ on $W$ with respective multiplicities $n_1,\ldots,n_{r+1}$, where $\ii = \sqrt{-1}$. Then $r \geq 1$. Set $d_j = \sum_{\ell = 1}^{j}n_\ell$ for $1 \leq j \leq r$. Every nonzero real tangent vector $X$ to $M$ satisfies
\begin{equation}\label{eq:delta-definition}
    H_g(X) \geq \frac{2}{S} > 0,
\end{equation}
where $S = \sum_{j = 1}^{r}(r - j + 1)(a_j - a_{j+1})d_j$.
\end{theorem}

\section{K\"ahler reduction and unitary flag manifolds}\label{sec:unitary-reduction}
We realize each invariant K\"ahler metric on a unitary partial flag manifold as a quotient metric. The defining equations of the level set then give an explicit lower bound for its holomorphic sectional curvature.

\subsection{Curvature of the quotient}
Let $(V,g_0,J_0,\omega_0)$ be a finite-dimensional Hermitian vector space with its flat K\"ahler structure. We normalize $g_0$ to be twice the real part of the Hermitian inner product and take $\omega_0(u,v) = g_0(J_0u,v)$. Let $K$ be a compact connected Lie group with Lie algebra $\mathfrak k$ acting linearly and unitarily on $V$. For $\eta \in \mathfrak k$, denote the infinitesimal action at $z \in V$ by $\eta_V(z) = \left.\frac{d}{dt}\right|_{t = 0}\exp(t\eta)\cdot z$. Let $\mu\colon V \to \mathfrak k^*$ be a $K$-equivariant moment map. For its component $\mu^\eta(z) = \langle\mu(z),\eta\rangle$, we use the convention $d\mu^\eta = \iota_{\eta_V}\omega_0$. Fix a regular value $c \in \mathfrak k^*$ invariant under the coadjoint action. Assume that $Z = \mu^{-1}(c)$ is nonempty and that $K$ acts freely on $Z$. Equivariance of $\mu$ ensures that $K$ preserves $Z$.
Write $\pi\colon Z \to Z/K$ for the quotient map. At $z \in Z$, the horizontal space $\mathcal H_z$ is the $g_0$-orthogonal complement of $T_z(K\cdot z)$ in $T_zZ$. Since $T_zZ = \ker d\mu_z$ and $d\mu^\eta_z(u) = g_0(J_0\eta_V(z),u)$, the normal space of $Z$ at $z$ is $J_0T_z(K\cdot z)$. Thus, $\mathcal H_z$ is the orthogonal complement in $T_zV$ of $T_z(K\cdot z) \oplus J_0T_z(K\cdot z)$. This direct sum is $J_0$-invariant, so its orthogonal complement is also $J_0$-invariant. The differential $d\pi_z$ identifies $\mathcal H_z$ with $T_{\pi(z)}(Z/K)$. For $u,v \in \mathcal H_z$, the reduced metric and complex structure satisfy
\begin{equation}\label{eq:quotient-structure}
    \begin{aligned}
        g_{\mathrm{red}}(d\pi_z u,d\pi_z v) & = g_0(u,v), \\
        J_{\mathrm{red}}d\pi_z u & = d\pi_z(J_0u).
    \end{aligned}
\end{equation}
Let $\II_Z$ be the second fundamental form of the real submanifold $Z \subset (V,g_0)$ and put $\|u\|_{g_0}^2 = g_0(u,u)$.

\begin{lemma}[Kobayashi, {\cite[Theorem~(1.3)]{Kobayashi}}]\label{lem:reduction}
Under the preceding assumptions, for every $z \in Z$ and every nonzero $u \in \mathcal H_z$,
\begin{equation}\label{eq:reduction}
    H_{g_{\mathrm{red}}}(d\pi_z u) = \frac{4\|\II_Z(u,u)\|_{g_0}^2}{\|u\|_{g_0}^4}.
\end{equation}
\end{lemma}
\begin{proof}
Here $J_0$ interchanges $T_z(K\cdot z)$ and the normal space of $Z$ while $d\pi_z$ is a complex-linear isometry on $\mathcal H_z$. These are the hypotheses of Kobayashi's submersion formula. Since $V$ is flat, that formula gives $H_{g_{\mathrm{red}}}(d\pi_z w) = 4\|\II_Z(w,w)\|_{g_0}^2$ for a unit horizontal vector $w$. Taking $w = u/\|u\|_{g_0}$ proves~\eqref{eq:reduction}. In terms of the moment map, the same identity follows from~\cite[Lemma~3.3, equations~(3.2)--(3.3)]{Zhang} applied to $\mu - c$, an equivariant moment map with zero set $Z$.
\end{proof}

When $Z$ lies on a sphere, the radial component of $\II_Z(u,u)$ alone gives a lower bound in~\eqref{eq:reduction}.

\begin{lemma}\label{lem:sphere}
For $\varrho > 0$, let $S_{\varrho} = \{z \in V \mid \|z\|_{g_0} = \varrho\}$. If $Z \subset S_{\varrho}$, then $H_{g_{\mathrm{red}}}(v) \geq 4/\varrho^2$ for every nonzero real tangent vector $v$ to $Z/K$.
\end{lemma}
\begin{proof}
Fix $z \in Z$ and $0 \ne u \in \mathcal H_z$. Since $T_zZ \subset T_zS_{\varrho} = (\mathbb R z)^\perp$, the radial vector $z$ is normal to $Z$. Choose a smooth curve $\gamma$ in $Z$ with $\gamma(0) = z$ and $\gamma'(0) = u$. Differentiating $g_0(\gamma(t),\gamma(t)) = \varrho^2$ twice at $t = 0$ gives $g_0(\gamma''(0),z) = -\|u\|_{g_0}^2$. Only the normal component $\II_Z(u,u)$ of $\gamma''(0)$ contributes to this pairing. Hence,
\begin{equation}\label{eq:sphere-pairing}
    g_0(\II_Z(u,u),z) = -\|u\|_{g_0}^2.
\end{equation}
The component of $\II_Z(u,u)$ along the unit radial vector $z/\varrho$ has length $\|u\|_{g_0}^2/\varrho$. Substituting this lower bound for $\|\II_Z(u,u)\|_{g_0}$ into~\eqref{eq:reduction} yields $H_{g_{\mathrm{red}}}(d\pi_z u) \geq 4/\varrho^2$. By~\eqref{eq:quotient-structure}, every nonzero tangent vector to $Z/K$ has a nonzero horizontal lift.
\end{proof}

\subsection{The partial flag manifold}
To apply Lemma~\ref{lem:sphere}, we describe the partial flag manifold as a quotient of a space of linear maps. Fix positive integers $r,N$ and dimensions satisfying
\begin{equation}\label{eq:dimensions}
    0 = d_0 < d_1 < \cdots < d_r < d_{r + 1} = N.
\end{equation}
Put $n_j = d_j - d_{j - 1}$ for $1 \leq j \leq r + 1$. Let $F(d_1,\ldots,d_r;\C^N)$ denote the manifold of flags $F_1 \subset \cdots \subset F_r \subset \C^N$ with $\dim_{\C} F_j = d_j$. Give each $\C^d$ its standard Hermitian inner product. Write $U(d)$ for its unitary group and $\mathfrak u(d)$ for the Lie algebra of skew-Hermitian matrices. Let $I_d$ be the identity on $\C^d$. For a linear map $A$ between Hermitian vector spaces, write $A^*$ for its adjoint and $\|A\|_{\HS}^2 = \Tr(A^*A)$ for its squared Hilbert--Schmidt norm. Set
\begin{equation}\label{eq:chain}
    \begin{aligned}
        V & = \bigoplus_{j = 1}^{r}\Hom(\C^{d_j},\C^{d_{j + 1}}), \\
        K & = \prod_{j = 1}^{r}U(d_j).
    \end{aligned}
\end{equation}
For $A = (A_1,\ldots,A_r) \in V$ and $X,Y \in T_AV \cong V$, the flat metric and K\"ahler form are
\begin{equation}\label{eq:chain-flat}
    \begin{aligned}
        g_0(X,Y) & = 2\operatorname{Re}\sum_{j = 1}^{r}\Tr(X_j^*Y_j), \\
        \omega_0(X,Y) & = 2\operatorname{Im}\sum_{j = 1}^{r}\Tr(X_j^*Y_j).
    \end{aligned}
\end{equation}
For $k = (k_1,\ldots,k_r) \in K$, put $k_{r + 1} = I_N$. The unitary action of $K$ on $V$ is given by
\begin{equation}\label{eq:chain-group-action}
    (k\cdot A)_j = k_{j + 1}A_j k_j^{-1}
\end{equation}
for $1 \leq j \leq r$. Let $A_0$ be the zero map from $\C^{d_0}$ to $\C^{d_1}$. For $\eta = (\eta_1,\ldots,\eta_r) \in \mathfrak k$, put $\eta_{r + 1} = 0$. Differentiating~\eqref{eq:chain-group-action} gives $(\eta_V(A))_j = \eta_{j + 1}A_j - A_j\eta_j$. Choose the moment map $\mu$ for this action so that $\mu(0) = 0$. Since $\eta_V(tA) = t\eta_V(A)$, the moment map identity gives $\frac{d}{dt}\mu^\eta(tA) = t\omega_0(\eta_V(A),A)$. Integration from $t = 0$ to $t = 1$ gives
\begin{equation}\label{eq:moment-sign}
    \mu^\eta(A)  = \frac12\omega_0(\eta_V(A),A) = -\ii\sum_{j = 1}^{r}\Tr\bigl((A_j^*A_j - A_{j - 1}A_{j - 1}^*)\eta_j\bigr).
\end{equation}
For any $c_1,\ldots,c_r > 0$, define the coadjoint-invariant functional $c \in \mathfrak k^*$ by
\[
    c(\eta) = -\ii\sum_{j = 1}^{r}c_j\Tr(\eta_j).
\]
By~\eqref{eq:moment-sign}, the defining equations of the level set $Z = \mu^{-1}(c)$ are
\begin{equation}\label{eq:level}
    A_j^*A_j - A_{j - 1}A_{j - 1}^* = c_j I_{d_j}
\end{equation}
for $1 \leq j \leq r$. Write $[A] = K\cdot A$ for the orbit of $A \in Z$. Let $U(N)$ act on $V$ by left multiplication on the last component,
\begin{equation}\label{eq:residual-action}
    u\cdot A = (A_1,\ldots,A_{r - 1},uA_r).
\end{equation}
This action preserves $Z$ since $(uA_r)^*(uA_r) = A_r^*A_r$. It commutes with~\eqref{eq:chain-group-action} and induces the action $u\cdot[A] = [u\cdot A]$ on $Z/K$. On the partial flag manifold, $U(N)$ acts by
\begin{equation}\label{eq:flag-group-action}
    u\cdot(F_1 \subset \cdots \subset F_r) = (uF_1 \subset \cdots \subset uF_r).
\end{equation}

\begin{proposition}\label{prop:unitary}
The quotient $Z/K$ is a smooth K\"ahler manifold. For $A \in V$ and $1 \leq j \leq r$, let $T_j(A) = A_r\cdots A_j\colon\C^{d_j} \to \C^N$. Write $\operatorname{im}T_j(A)$ for the subspace $T_j(A)(\C^{d_j})$ of $\C^N$. The map
\begin{equation}\label{eq:quotient-flag-map}
    \begin{aligned}
        \phi\colon Z/K &\longrightarrow F(d_1,\ldots,d_r;\C^N), \\
        [A] &\longmapsto \bigl(\operatorname{im}T_1(A) \subset \cdots \subset \operatorname{im}T_r(A) \subset \C^N\bigr)
    \end{aligned}
\end{equation}
is a $U(N)$-equivariant biholomorphism for the actions in~\eqref{eq:residual-action} and~\eqref{eq:flag-group-action}.
\end{proposition}
\begin{proof}
The level set $Z$ is nonempty. Indeed, let $\iota_j\colon\C^{d_j} \to \C^{d_{j+1}}$ be the coordinate inclusion. Starting with $A_0 = 0$, define $A_j = \iota_j\circ(A_{j-1}A_{j-1}^* + c_j I_{d_j})^{1/2}$ successively using the positive square root. Since $\iota_j^*\iota_j = I_{d_j}$, these maps satisfy~\eqref{eq:level}. For an arbitrary $A \in Z$ and $v \in \C^{d_j}$, the same equation gives
\[
    \|A_jv\|^2 = \|A_{j-1}^*v\|^2 + c_j\|v\|^2.
\]
The right-hand side is positive for $v \ne 0$, so every $A_j$ is injective. If $k \in K$ fixes $A$, then $A_rk_r^{-1} = A_r$ implies $k_r = I_{d_r}$. Descending induction in~\eqref{eq:chain-group-action} gives $k_j = I_{d_j}$ for all $j$. Thus, $K$ acts freely on $Z$. If $\eta \in \mathfrak k$ annihilates the image of $d\mu_A$, then $g_0(J_0\eta_V(A),X) = 0$ for all $X \in T_AV$. It follows that $\eta_V(A) = 0$, so $\eta$ lies in the Lie algebra of the trivial stabilizer. Hence, $d\mu_A$ is surjective at every $A \in Z$. Reduction at this regular value by the free action of the compact group $K$ gives a smooth K\"ahler quotient.
Let $\Omega \subset V$ be the open set on which every $A_j$ is injective. For $A \in \Omega$, the composite $T_j(A)$ is injective and the identity $T_j(A) = T_{j+1}(A)A_j$ for $j < r$ gives $\operatorname{im}T_j(A) \subset \operatorname{im}T_{j+1}(A)$. Thus, these images form a flag of the prescribed dimensions. Let $\widetilde\phi\colon\Omega \to F(d_1,\ldots,d_r;\C^N)$ send $A$ to this flag. The Pl\"ucker coordinates of its $j$th subspace are the maximal minors of $T_j(A)$. These minors are polynomial in the entries of $A_1,\ldots,A_r$. At each $A \in \Omega$, at least one maximal minor of $T_j(A)$ is nonzero since $T_j(A)$ has rank $d_j$. In the corresponding Pl\"ucker chart, the coordinates are ratios of these minors. Hence, $\widetilde\phi$ is holomorphic. We have already shown that $Z \subset \Omega$. Since $T_j(k\cdot A) = T_j(A)k_j^{-1}$, the restriction of $\widetilde\phi$ to $Z$ is constant on $K$-orbits and descends to the map $\phi$ in~\eqref{eq:quotient-flag-map}. The quotient description of Donagi and Sharpe~\cite[Sections~2.1--2.2]{DS} gives bijectivity for $c_1,\ldots,c_r > 0$, with the flag identified by the Pl\"ucker coordinates above. The identity $T_j(u\cdot A) = uT_j(A)$ gives $U(N)$-equivariance. Since $U(N)$ acts transitively on the flag manifold, bijectivity and equivariance imply that it also acts transitively on $Z/K$ and that corresponding points have the same stabilizer. Thus, $\phi$ is a diffeomorphism. Finally, differentiating $\widetilde\phi|_Z = \phi\circ\pi$ along the horizontal space $\mathcal H_A$ gives
\[
    d\widetilde\phi_A|_{\mathcal H_A} = d\phi_{[A]}\circ(d\pi_A|_{\mathcal H_A}).
\]
The left-hand side is complex-linear because $\widetilde\phi$ is holomorphic. By~\eqref{eq:quotient-structure}, $d\pi_A|_{\mathcal H_A}$ is a complex-linear isomorphism. Hence, $d\phi_{[A]}$ is complex-linear at every $[A] \in Z/K$, so $\phi$ is a biholomorphism.
\end{proof}

\subsection{The metric and curvature bound}
Proposition~\ref{prop:unitary} identifies the quotient as a complex manifold. We first bound the curvature of its reduced metric and then determine which level set induces a prescribed invariant K\"ahler metric. For the constants in~\eqref{eq:level}, put
\begin{equation}\label{eq:S-definition}
    S = \sum_{j = 1}^{r}(r - j + 1)c_j d_j.
\end{equation}
\begin{lemma}\label{lem:unitary-metric}
For every choice of $c_1,\ldots,c_r > 0$, the reduced metric on $Z/K$ satisfies
\begin{equation}\label{eq:unitary-bound}
    H_{g_{\mathrm{red}}}(v) \geq \frac{2}{S} > 0
\end{equation}
for every nonzero real tangent vector $v$ to $Z/K$. For every $U(N)$-invariant K\"ahler metric $g$ on $F(d_1,\ldots,d_r;\C^N)$, there is a choice of these parameters for which the corresponding map $\phi$ satisfies $\phi^*g = g_{\mathrm{red}}$.
\end{lemma}
\begin{proof}
Taking traces in~\eqref{eq:level} gives $\|A_j\|_{\HS}^2 - \|A_{j-1}\|_{\HS}^2 = c_jd_j$ for every $A \in Z$. Since $A_0 = 0$, summing first up to $j$ and then over $j$ gives
\begin{equation}\label{eq:trace}
    g_0(A,A) = 2\sum_{j = 1}^{r}\|A_j\|_{\HS}^2 = 2\sum_{j = 1}^{r}\sum_{\ell = 1}^{j}c_\ell d_\ell = 2S.
\end{equation}
Thus, $Z$ lies on the sphere of radius $\sqrt{2S}$. Lemma~\ref{lem:sphere} gives $H_{g_{\mathrm{red}}} \geq 4/(2S) = 2/S$. Now let $\omega$ be the K\"ahler form of a given $U(N)$-invariant K\"ahler metric $g$. Let $F_\bullet$ be the standard coordinate flag and let $p(u) = u\cdot F_\bullet$ be its orbit map. Write $E_j$ for the span of the coordinate vectors numbered $d_{j-1}+1$ through $d_j$. Since $\omega$ is $SU(N)$-invariant, the type $A$ parametrization of invariant K\"ahler forms~\cite[Section~2.5]{AD} gives real numbers $a_1 > \cdots > a_{r+1}$ such that the Hermitian endomorphism $\Lambda$ acting by $a_j$ on $E_j$ satisfies
\begin{equation}\label{eq:unitary-flag-form}
    (p^*\omega)_{I_N}(X,Y) = \ii\Tr(\Lambda[X,Y]).
\end{equation}
This formula holds for all $X,Y \in \mathfrak u(N)$ because scalar matrices act trivially on flags and their Lie brackets vanish. Adding a scalar multiple of $I_N$ to $\Lambda$ does not change the form since $\Tr([X,Y]) = 0$. We may therefore take $a_{r+1} = 0$ and choose $c_j = a_j - a_{j+1}$.
To compare the two forms at $F_\bullet$, identify $\C^{d_j}$ with $E_1 \oplus \cdots \oplus E_j$. For $1 \leq j \leq r$, define $A_j\colon\C^{d_j} \to \C^{d_{j+1}}$ by $A_jv = (a_\ell - a_{j+1})^{1/2}v$ on each $E_\ell$ with $\ell \leq j$. On $E_\ell$, the operators $A_j^*A_j$ and $A_{j-1}A_{j-1}^*$ act by $a_\ell - a_{j+1}$ and $a_\ell - a_j$, respectively. Their difference is $c_j I_{d_j}$, so $A \in Z$. The images of its composites form the standard flag and the formula for $A_r$ gives
\begin{equation}\label{eq:spectrum}
    \begin{aligned}
        \phi([A]) & = F_\bullet, \\
        A_rA_r^* & = \Lambda = \diag(a_1 I_{n_1},\ldots,a_{r+1} I_{n_{r+1}}).
    \end{aligned}
\end{equation}
Let $q_A(u) = [u\cdot A]$ be the orbit map into $Z/K$. Its lift $u \mapsto u\cdot A$ has differential $X \mapsto (0,\ldots,0,XA_r)$ at $I_N$. The pullback of $\omega_{\mathrm{red}}$ to $Z$ is the restriction of $\omega_0$. Thus, equations~\eqref{eq:chain-flat} and~\eqref{eq:spectrum} give
\begin{equation}\label{eq:residual-form-calculation}
    \begin{aligned}
        (q_A^*\omega_{\mathrm{red}})_{I_N}(X,Y)
        & = 2\operatorname{Im}\Tr\bigl((XA_r)^*YA_r\bigr) \\
        & = \ii\Tr\bigl(A_rA_r^*[X,Y]\bigr) = (p^*\omega)_{I_N}(X,Y).
    \end{aligned}
\end{equation}
By equivariance of $\phi$, we have $p = \phi\circ q_A$. Since the action is transitive, $(dq_A)_{I_N}$ is surjective. Hence,~\eqref{eq:residual-form-calculation} gives $\omega_{\mathrm{red}} = \phi^*\omega$ at $[A]$. The $U(N)$-action preserves $Z$ and $\omega_0$ and commutes with $K$, so it preserves $\omega_{\mathrm{red}}$. Both forms are therefore $U(N)$-invariant and agree everywhere by transitivity. Since $\phi$ is holomorphic by Proposition~\ref{prop:unitary}, their metrics satisfy $\phi^*g = g_{\mathrm{red}}$.
\end{proof}

\section{Totally geodesic embeddings of flag manifolds}\label{sec:totally-geodesic}
We first identify $(M,J,g)$ with the adjoint orbit through the element $\xi_\omega$ determined by~\eqref{eq:representative}. Then, we show that the inclusion of two adjoint orbits through the same element is holomorphic, isometric and totally geodesic when both orbit metrics are defined using the same invariant inner product.

\subsection{Invariant K\"ahler forms}
For $\xi \in \g$, let $o = eC_G(\xi)$ be the identity coset in $G/C_G(\xi)$. Set $L = \ad(\xi)$ and let $\m = (\ker L)^{\perp_Q}$. The differential of the quotient map identifies $\m$ with $T_o(G/C_G(\xi))$. Since $L$ is $Q$-skew-adjoint and injective on $\m$, the operator $-L^2$ is positive definite there. Denote its positive square root on $\m$ by $D$. Let $(J_\xi,g_\xi,\omega_\xi)$ be the canonical K\"ahler structure on the orbit. With the sign convention of~\eqref{eq:ambient-representative}, the Kirillov--Kostant--Souriau formula~\cite[Chapter~5, p.~29]{DHJH} and the description of the complex structure in~\cite[Section~2.2.2]{DonaldsonLie} give
\begin{equation}\label{eq:orbit-data}
    \begin{aligned}
        (\omega_\xi)_o(X,Y) & = -Q(LX,Y), \\
        (J_\xi)_o & = -LD^{-1}, \\
        (g_\xi)_o(X,Y) & = Q(DX,Y)
    \end{aligned}
\end{equation}
for $X,Y \in \m$.

\begin{lemma}\label{lem:invariant-form}
For the fixed inner product $Q$, the K\"ahler form $\omega$ of $(M,J,g)$ determines a unique element $\xi_\omega \in \g$ satisfying~\eqref{eq:representative}. Its centralizer in $G$ is $H$. The map
\[
    aH \longmapsto \Ad_a\xi_\omega
\]
is a $G$-equivariant holomorphic isometry from $(M,J,g)$ onto $(\Ad(G)\xi_\omega,J_{\xi_\omega},g_{\xi_\omega})$, where $J_{\xi_\omega}$ and $g_{\xi_\omega}$ are defined by~\eqref{eq:orbit-data}.
\end{lemma}
\begin{proof}
By the Kirillov--Kostant--Souriau theorem and its corollary in~\cite[Section~II]{Alekseevsky97}, there is a functional $\ell_\omega \in \g^*$ whose coadjoint stabilizer is $H$ and which satisfies $(p_G^*\omega)_e(X,Y) = -\ell_\omega([X,Y])$ for all $X,Y \in \g$. The functional is unique because $[\g,\g] = \g$. Its $Q$-dual $\xi_\omega$ satisfies~\eqref{eq:representative}. Since $Q$ identifies the adjoint and coadjoint actions, $C_G(\xi_\omega) = H$. Thus, the map $\Psi(aH) = \Ad_a\xi_\omega$ is a $G$-equivariant diffeomorphism onto $\Ad(G)\xi_\omega$. Equations~\eqref{eq:representative} and~\eqref{eq:orbit-data} give $\Psi^*\omega_{\xi_\omega} = \omega$ at $eH$ and hence everywhere by $G$-invariance. By~\cite[Section~IV.5, corollary and subsequent remark]{Alekseevsky97}, an invariant symplectic form is K\"ahler precisely for the invariant complex structure corresponding to the chamber containing its parameter. Since $J$ and $\Psi^*J_{\xi_\omega}$ both make $\omega$ K\"ahler, they coincide. Therefore, $\Psi$ is holomorphic and $\Psi^*g_{\xi_\omega} = \omega(\cdot,J\cdot) = g$.
\end{proof}

\subsection{Totally geodesic embeddings of adjoint orbits}
Let $G \subset \widetilde G$ be compact connected Lie groups and let $Q$ be an $\Ad(\widetilde G)$-invariant inner product on $\widetilde{\g}$. For $\xi \in \g$, the formulas in~\eqref{eq:orbit-data} show that the natural inclusion of $G/C_G(\xi)$ into $\widetilde G/C_{\widetilde G}(\xi)$ is holomorphic and isometric for the K\"ahler structures defined by $Q|_{\g}$ and $Q$, respectively. To prove total geodesy, we use the following criterion of Alekseevsky and Nikonorov.

\begin{lemma}[Alekseevsky--Nikonorov, {\cite[Section~3, Corollary~2]{AN}}]\label{lem:subgroup-orbit}
Let $\widetilde G$ be a compact connected Lie group and $\widetilde H$ a closed subgroup. Equip $\widetilde G/\widetilde H$ with a $\widetilde G$-invariant Riemannian metric and put $o = e\widetilde H$. Write $\widetilde{\g}$ and $\widetilde{\h}$ for the Lie algebras of $\widetilde G$ and $\widetilde H$. Let $\widetilde{\m}$ be the orthogonal complement of $\widetilde{\h}$ with respect to an $\Ad(\widetilde G)$-invariant inner product on $\widetilde{\g}$ and identify $\widetilde{\m}$ with $T_o(\widetilde G/\widetilde H)$. For a compact connected subgroup $G \subset \widetilde G$ with Lie algebra $\g$, set $\m = \g \cap \widetilde{\m}$. Suppose that
\[
    \g = (\g \cap \widetilde{\h}) \oplus \m.
\]
Let $\mathfrak n$ be the orthogonal complement of $\m$ in $\widetilde{\m}$ for the Riemannian metric at $o$. Denote by $P_{\m}$ the projection onto $\m$ along $\widetilde{\h} \oplus \mathfrak n$. The orbit $G\cdot o$ is totally geodesic if and only if the endomorphism $X \mapsto P_{\m}[Z,X]$ of $\m$ is skew-adjoint for the metric at $o$ for every $Z \in \mathfrak n$.
\end{lemma}

\begin{proposition}\label{prop:totally-geodesic}
Let $G \subset \widetilde G$ be compact connected Lie groups with Lie algebras $\g \subset \widetilde{\g}$. Fix $\xi \in \g$ and an $\Ad(\widetilde G)$-invariant inner product $Q$ on $\widetilde{\g}$. Equip $G/C_G(\xi)$ and $\widetilde G/C_{\widetilde G}(\xi)$ with the K\"ahler structures in~\eqref{eq:orbit-data} determined by $\xi$ using $Q|_{\g}$ and $Q$, respectively. Then the map
\begin{equation}\label{eq:inclusion}
    \begin{aligned}
        \Phi \colon G/C_G(\xi) &\longrightarrow \widetilde G/C_{\widetilde G}(\xi), \\
        aC_G(\xi) &\longmapsto aC_{\widetilde G}(\xi)
    \end{aligned}
\end{equation}
is a holomorphic, isometric, totally geodesic embedding.
\end{proposition}
\begin{proof}
Let $o = eC_G(\xi)$ and $\widetilde o = \Phi(o)$. Under the adjoint orbit identifications, $\Phi$ is the inclusion $\Ad(G)\xi \hookrightarrow \Ad(\widetilde G)\xi$. The orbit $\Ad(G)\xi$ is embedded because $G$ is compact. To compare the K\"ahler structures at $o$ and $\widetilde o$, set $\widetilde L = \ad(\xi)$ on $\widetilde{\g}$ and write $\widetilde{\h} = \ker\widetilde L$ and $\widetilde{\m} = \widetilde{\h}^{\perp_Q}$. For $X,Y \in \g$ and $Z \in \g^{\perp_Q}$, invariance of $Q$ gives
\begin{equation}\label{eq:complement}
    Q([X,Z],Y) = -Q(Z,[X,Y]) = 0.
\end{equation}
Thus, $[\g,\g^{\perp_Q}] \subset \g^{\perp_Q}$. Since $\xi \in \g$, the operator $\widetilde L$ preserves both summands of $\widetilde{\g} = \g \oplus \g^{\perp_Q}$. Its kernel therefore splits along these summands, as does its $Q$-orthogonal complement. Setting $\m = \g \cap \widetilde{\m}$ and $\mathfrak n = \g^{\perp_Q} \cap \widetilde{\m}$ gives the $Q$-orthogonal decompositions
\begin{equation}\label{eq:orthogonal-split}
    \begin{aligned}
        \g & = (\g \cap \widetilde{\h}) \oplus \m, \\
        \widetilde{\m} & = \m \oplus \mathfrak n.
    \end{aligned}
\end{equation}
Under the tangent space identifications at $o$ and $\widetilde o$, the differential $d\Phi_o$ is the inclusion $\m \hookrightarrow \widetilde{\m}$. Write $L = \widetilde L|_{\g}$ and denote the positive square roots of $-L^2$ on $\m$ and $-\widetilde L^2$ on $\widetilde{\m}$ by $D$ and $\widetilde D$, respectively. Since $\widetilde L$ preserves the orthogonal summands $\m$ and $\mathfrak n$, so does $\widetilde D$. By uniqueness of the positive square root, $\widetilde D|_{\m} = D$. Thus, the ambient complex structure restricts to $-LD^{-1}$ on $\m$ and the ambient metric restricts to $Q(D\cdot,\cdot)$ there. These are precisely the complex structure and metric in~\eqref{eq:orbit-data} for $G/C_G(\xi)$. By $G$-equivariance, $\Phi$ is holomorphic and isometric everywhere.
It remains to verify the criterion in Lemma~\ref{lem:subgroup-orbit}. The first equality in~\eqref{eq:orthogonal-split} gives the required decomposition of $\g$. Let $\widetilde g$ denote the metric on $\widetilde G/C_{\widetilde G}(\xi)$. The second equality is also orthogonal for $\widetilde g_{\widetilde o}$ because $\widetilde g_{\widetilde o}(X,Z) = Q(\widetilde D X,Z) = 0$ for $X \in \m$ and $Z \in \mathfrak n$. Hence, $\mathfrak n$ is the normal space to $G\cdot\widetilde o$ at $\widetilde o$. The projection $P_{\m}$ along $\widetilde{\h} \oplus \mathfrak n$ is the $Q$-orthogonal projection onto $\m$. For $Z \in \mathfrak n$ and $X \in \m$, equation~\eqref{eq:complement} gives $[Z,X] \in \g^{\perp_Q}$. Therefore, $P_{\m}[Z,X] = 0$. The zero endomorphism is skew-adjoint, so Lemma~\ref{lem:subgroup-orbit} proves that $\Phi$ is totally geodesic.
\end{proof}

\section{Curvature of generalized flag manifolds}\label{sec:general-curvature}

\begin{proof}[Proof of Theorem~\ref{thm:main}]
Use the trace inner product $Q(X,Y) = -\Tr(XY)$ on $\mathfrak u(W)$ and its restriction to $\g$. By Lemma~\ref{lem:invariant-form}, $H = C_G(\xi_\omega)$ and the orbit metric on $G/C_G(\xi_\omega)$ agrees with $g$. Applying Proposition~\ref{prop:totally-geodesic} with $\widetilde G = U(W)$ and $\xi = \xi_\omega$ gives the holomorphic, isometric and totally geodesic embedding
\begin{equation}\label{eq:general-embedding}
    \begin{aligned}
        \Phi\colon (M,g) &\longrightarrow (M_W,g_W), \\
        aH &\longmapsto aC_{U(W)}(\xi_\omega).
    \end{aligned}
\end{equation}
For the curvature estimate, $C_G(\xi_\omega) = H \ne G$, so $\xi_\omega$ is not scalar. Let $E_j$ be the eigenspace of $-\ii\xi_\omega$ for $a_j$ and put $F_j = E_1 \oplus \cdots \oplus E_j$. The map $uC_{U(W)}(\xi_\omega) \mapsto (uF_1 \subset \cdots \subset uF_r)$ identifies $M_W$ with the partial flag manifold. In the order $a_1 > \cdots > a_{r+1}$, the strictly lower blocks give holomorphic tangent coordinates on the flag manifold. By~\eqref{eq:orbit-data}, $J_W$ multiplies these blocks by $\ii$, so the identification is holomorphic. Comparing~\eqref{eq:ambient-representative} with~\eqref{eq:unitary-flag-form} gives $\Lambda = -\ii\xi_\omega$ up to addition of a scalar multiple of $I_W$. Thus, the proof of Lemma~\ref{lem:unitary-metric} gives $c_j = a_j - a_{j+1}$ for $g_W$. With these parameters, the lemma gives $H_{g_W} \geq 2/S$ with $S$ as in Theorem~\ref{thm:main}. To transfer the estimate to $M$, write $R_W$ for the curvature tensor of $g_W$ and $\II_\Phi$ for the second fundamental form of $\Phi$. For tangent vectors $X_j$ at a point of $M$, set $U_j = d\Phi(X_j)$. With convention~\eqref{eq:Riemannian-curvature-tensor}, the Gauss equation is
\begin{equation}\label{eq:gauss-embedding}
    \begin{aligned}
        g(R(X_1,X_2)X_3,X_4)
        & = g_W(R_W(U_1,U_2)U_3,U_4) \\
        &\phantom{{}={}} + g_W(\II_\Phi(X_1,X_4),\II_\Phi(X_2,X_3)) \\
        &\phantom{{}={}} - g_W(\II_\Phi(X_1,X_3),\II_\Phi(X_2,X_4)).
    \end{aligned}
\end{equation}
Total geodesy gives $\II_\Phi = 0$. For any nonzero real tangent vector $X$ to $M$, substitute $(X,JX,JX,X)$ in~\eqref{eq:gauss-embedding}. Since $d\Phi(JX) = J_W d\Phi(X)$ and $\Phi$ preserves lengths, division by $g(X,X)^2$ gives
\begin{equation}\label{eq:explicit-bound}
    H_g(X) = H_{g_W}(d\Phi(X)) \geq \frac{2}{S} > 0.
\end{equation}
\end{proof}

\section*{Disclosure of AI assistance}
The author used OpenAI's ChatGPT to explore the mathematics, check computations and revise the manuscript.

\end{document}